\documentclass{article}

\usepackage{amsmath,amssymb,mathrsfs,theorem,subfiles}
\usepackage{graphicx}
\usepackage{subcaption}
\usepackage{tabu}
\usepackage[table]{xcolor}
\usepackage{diagbox}
\usepackage{multirow}
\usepackage{url}
\usepackage{wrapfig}

\newtheorem{definition}{Definition}[section]

\newtheorem{rem}[definition]{Remark}

\newcommand{\X}{\mathcal{X}}
\newcommand{\C}{\mathcal{C}}

\usepackage[top=30truemm, bottom=30truemm,left=30truemm,right=30truemm]{geometry}

\begin{document}

\title{Analysis on MathSciNet database: some preliminary results}

\author{Serge Richard\footnote{Supported by the grant\emph{Topological invariants
through scattering theory and noncommutative geometry} from Nagoya University,
and by JSPS Grant-in-Aid for scientific research C no 18K03328, and on
leave of absence from Univ.~Lyon, Universit\'e Claude Bernard Lyon 1, CNRS UMR 5208,
Institut Camille Jordan, 43 blvd.~du 11 novembre 1918, F-69622 Villeurbanne cedex,
France.},
Qiwen Sun}

\date{\small}
\maketitle
\vspace{-1cm}

\begin{quote}
\emph{
\begin{enumerate}
\item[] Graduate school of mathematics, Nagoya University,
Chikusa-ku, Nagoya 464-8602, Japan
\item[]\emph{E-mails:} richard@math.nagoya-u.ac.jp, m17046h@math.nagoya-u.ac.jp
\end{enumerate}
}
\end{quote}

\begin{abstract}
In this paper we initiate some investigations
on MathSciNet database. For many mathematicians
this website is used on a regular basis, but surprisingly
except for the information provided by MathSciNet itself,
there exist almost no independent investigations
or independent statistics on this database.
This current research has been triggered by a rumor:
do international collaborations increase the number of citations
of an academic work in mathematics\;\!?
We use MathSciNet for providing some information
about this rumor, and more generally pave the way for further
investigations on or with MathSciNet.
\end{abstract}


\smallskip



\section{Introduction}

This research has been triggered by a rumor: do international
collaborations increase the number of citations in mathematics\;\!?
By looking at the existing studies about this question in various
fields of research \cite{LFSEM,SWBA}, the easy and naive answer would be positive.
For example, in reference \cite{SWBA} mathematics is one of the eight disciplines
considered, but the database used (Scopus) is not specific to mathematics.
Therefore, is it possible to get a more precise and deeper answer
by looking at a database specific to mathematics\;\!?
Among the very few choices, we decided to concentrate
on one database which is familiar to most mathematicians: MathSciNet.
In addition, since this database contains a lot of information,
the initial question has been generalized to a more natural one:
is it possible to identify some predictors at the time of publication of a paper
which will certainly impact the future citations of this paper\;\!?

Let us state it immediately, the information we provide and our results
are very partial and have to be considered as preliminary.
The reason is the following: MathSciNet is equipped with several tools for
searching and arranging publications, but it does not provide any tool for extracting
data from the website. Also, MathSciNet does not allow any automated searching or downloading.
As a consequence, even though the database contains information on 3.6 million works,
the total number of items that we can fully analyze is very limited.
Because extracting
and labelling information from the database is difficult and time consuming,
we had to design our experiments very carefully and modestly. However,
we hope that this initial work will justify a better and easier access to the database
in the near future. For completeness, let us mention that about 20 years ago
MathSciNet has been extensively used for the elaboration of some graphs about mathematics and mathematicians \cite{G1,G2}, but the authors are not aware of more recent works.

Taking these very strong restrictions into account, let us briefly describe our approach.
Since our initial motivation was to study the effect of international collaborations,
we firstly fix one target country, which is going to be Japan.
Also, since the response is related to the number of citations,
it obviously depends on time: A paper published 20 years ago and a paper published
last year will never have the same amount of time for being cited.
To eliminate this effect, we have decided to consider only papers published
in the same year, which we fixed to be 2009.
As a consequence, all predictors considered are from 2009.
In this restricted setting, the total number of publications 
available on MathSciNet is 3907 (in December 2018).

In reference \cite{SWBA} the predictors employed for each paper were the number of authors,  the number of countries (given by the affiliations of the authors), the number of references, the year of publication, and the impact factor of the journal.
Additional specific information (as for example economic and demographic data) related to the countries of the authors were also used in their analysis.
For our investigations we stuck to the information provided by MathSciNet,
but tried to exploit all of them. For example, we also took into account
the academic age of the authors, and the Primary Mathematics Subject Classification
(and implemented this information in a simpler classification scheme).
In the same vein, we reported the number of pages of each paper, or the existence
(or non-existence) of a review text in MathSciNet.
At the end of the day we have identified 13 predictors available for each paper
at the time of publication.

Even if our initial number of target papers was only 3907, it was beyond our means to
collect manually these 13 predictors, the response (the number of citations),
and some additional information which were finally not exploited.
For that reason we implemented a stratified sampling and collected the
full set of information on about 300 documents. This number is certainly too small
to get any firm answer, but nevertheless it allows to draw a preliminary
picture of the relations between these predictors and the response.

Let us now be more specific about the content of the paper.
In Section \ref{back} we provide the necessary information
about MathSciNet and present a few general features about
the publications for the last 20 years. Some of these features
motivate our subsequent choice of some predictors.
Our sample for the data and the list of predictors with some explanations
are also introduced in this section.
Based on our sample we gather in Section \ref{pre} some preliminary results
and establish some relations between the predictors and the response.
Note that in Figure \ref{fig:bx1} a relation between international
collaborations and the number of citations is visible. However,
the importance of this correlation will be weakened later on when
all predictors will be put on an equal footing.

The content of Section \ref{invest} corresponds to the core of our investigations.
Because of the diversity of the predictors we opted for an approach based on
tree-based methods \cite{BFOS}. Indeed, unlike the approach provided in \cite{SWBA}
we do not want to consider some linear relations between the predictors and the
response, but prefer an approach which divides the parameter space into several
regions and associates to each region a local response. Clearly, due to the small
size of our sample, the local response can not be the exact number of citations
but only an approximation of this number. For that reason, we have divided the response
into 3 classes according to the number of citations: the \emph{low level} for 0 or 1
citation, the \emph{median level} for a number of citations between 2 and 10,
and the \emph{high level} for a higher number of citations.
With all these information, a tree classifier has been established by using
the R programming. The outcomes are provided in this section.
We also discuss the relative importance of the predictors, and provide a
ranking of them according to two criteria. Let us immediately mention it: the predictor corresponding to the number of countries involved in the research project comes last,
for both rankings. This is one of the surprising results of our investigations.
Finally, in the last section we provide a conclusion for our investigations.

As already emphasized above, our conclusions are only preliminary, mainly
because of the lack of an easy access to the database. Nevertheless, we believe
that the approach should be implemented on a larger scale, and that
some of the outcomes deserve further investigations.
From the authors' point of view, MathSciNet is an incredibly rich and powerful
database which has been underinvestigated. We hope that our investigations
will pave the way for future research.

\section{The MathSciNet database}\label{back}

In this section we provide a few information about MathSciNet.
Some of the features presented motivate the choice of future predictors.
We also introduce our sample, and the list of predictors with some explanations.

\subsection{Some preliminary features}\label{sec:database}

\begin{wrapfigure}{R}{0.5\textwidth}
\centering
\includegraphics[scale = 0.34]{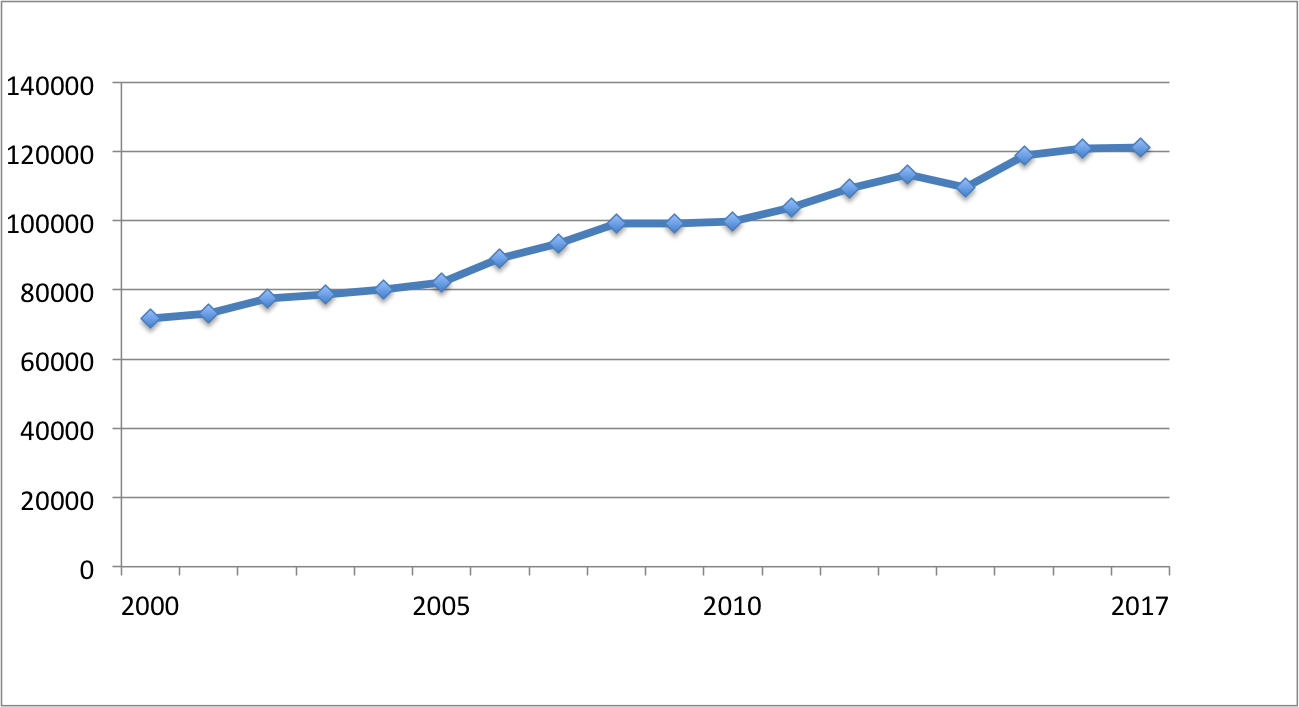}
\caption{Publications per year}
\label{fig:glo}
\end{wrapfigure}

This research is based on data about mathematics which are extracted
from the most well-known database in this field, which is \emph{MathSciNet}.
MathSciNet is an electronic database operated by the American Mathematical Society.
It indexes more than 3.6 million works in mathematics. 
Each year, more than 100,000 new items are added
to the database.
MathSciNet is equipped with several tools for searching and arranging publications.
One can easily search the works by a lot of search terms,
for example the MR number, the Mathematics Subject Classification (MSC),
the name of the author(s) or the name of the reviewer.
One can also fix a period of publication and list all the items published in the period.

\begin{wrapfigure}{L}{0.5\textwidth}
\centering
\includegraphics[scale = 0.35]{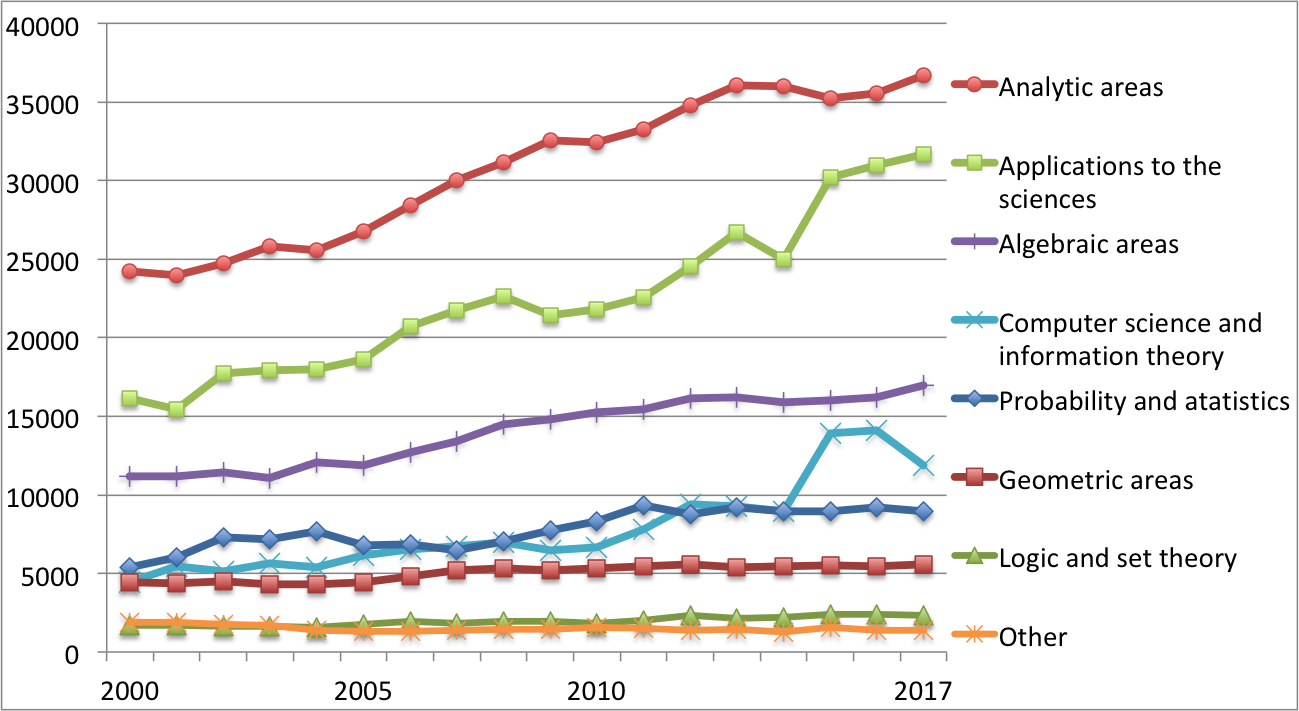}
 \caption{Publications by subject}
 \label{fig:pri}
\end{wrapfigure}

As a warm-up, we exhibit a few general features of the data
contained in MathSciNet. Note that all the subsequent pictures and tables have been
prepared by the authors based on data extracted from the database.
First of all, let us have a look at the global variation of publications
in mathematics for the recent 18 years.
From Figure \ref{fig:glo} one observes that the number of items collected
by the database is constantly increasing. This certainly reflects the general trend
about the number of publications in mathematics for this period of time.
However, if we take a look at Figure \ref{fig:pri} one observes that
this evolution highly depends on the subjects. Note that for this figure we have
used a combined subject classification, see Remark \ref{rem_MSC}.
Since this variation might have an impact on the number of citations,
we shall use the subject as a predictor in our subsequent analysis.

\begin{rem}\label{rem_MSC}
The Mathematics Subject Classification (MSC)  is an alphanumerical classification scheme
extensively used in mathematics. In MathSciNet, each item receives precisely one primary MSC (usually provided by the authors) which describes the principal subject of the paper. The current MSC2010 contains more than 60 subjects, which is too big compared
to the size of the sample we shall use subsequently.
For that reason, we shall rely on a simpler classification scheme provided by \cite{RD},
which classifies all items in 8 distinct subjects: logic and set theory, algebraic areas, geometric areas, analytic areas, probability and statistics,
computer science and information theory, applications to the sciences, others.
The primary MSC can be distributed into these 8 subjects based on a list provided
by \cite{RD}.
\end{rem}

\begin{figure}
\centering
\begin{minipage}[b]{0.4\textwidth}
\includegraphics[scale = 0.33]{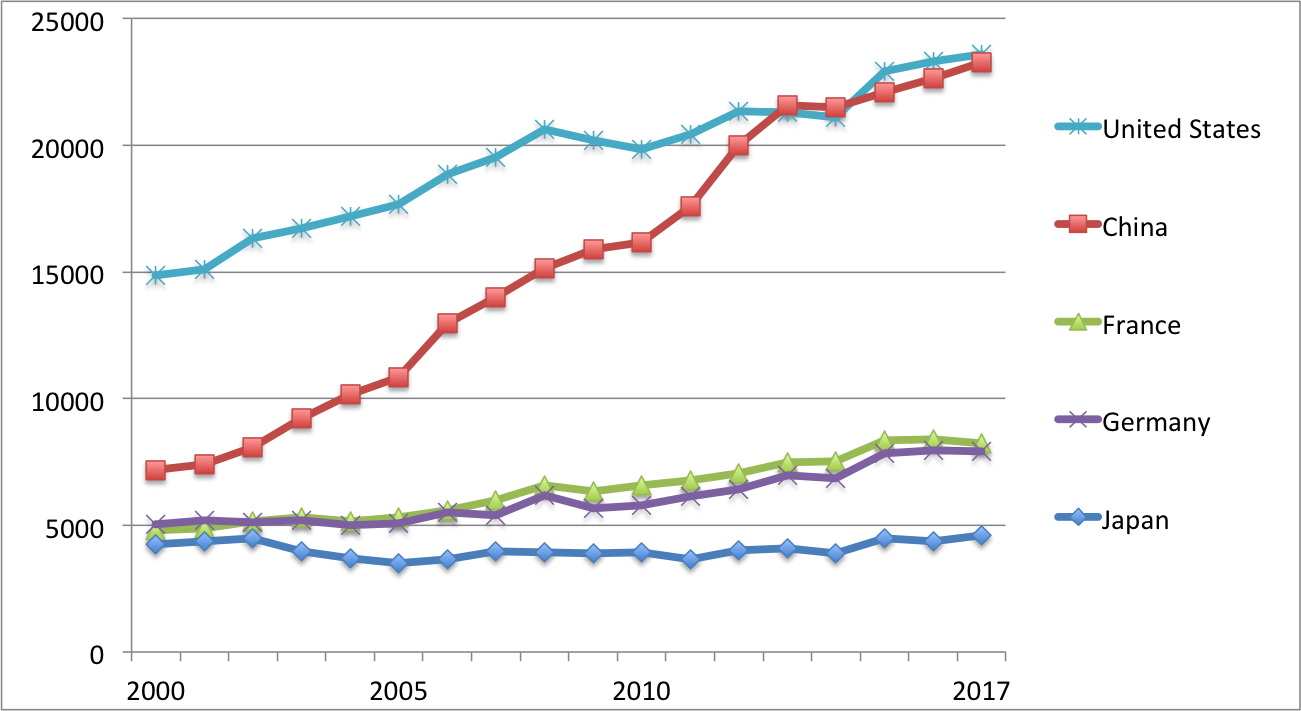} 
\caption{Publications with at least one author from a given country}
\label{fig:con}
\end{minipage}%
\hfill
\begin{minipage}[b]{0.5\textwidth}
\includegraphics[scale = 0.33]{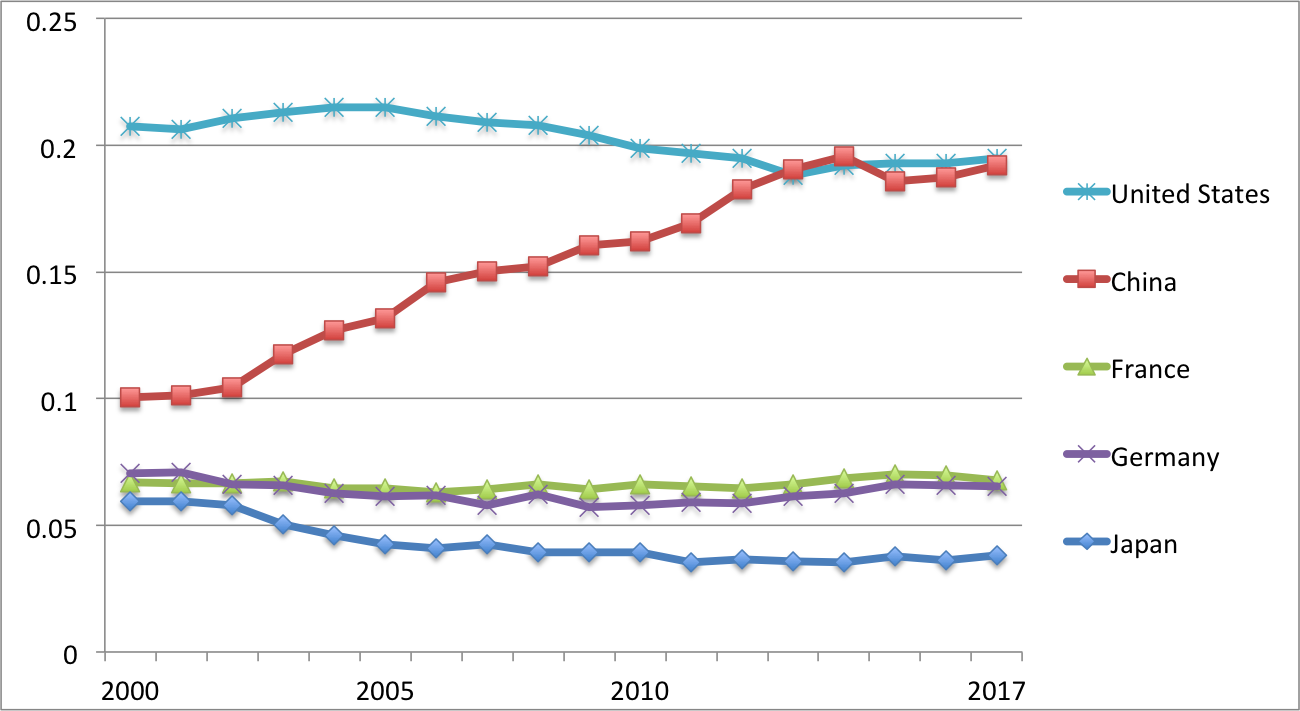}
\caption{Relative publications with at least one author from a given country}
\label{fig:rcon}
\end{minipage}
\end{figure}

Let us come to a more local criterion. For that purpose we shall use
the affiliation (university, research institute, etc.) provided by each author
for each publication.  Figure \ref{fig:con} shows the variation
of the numbers of publications with at least one author from an institute
in one of the following countries: Japan, China, France, Germany, and the United States.
These numbers are increasing except for Japan which is more of less stable.
On the other hand, if we divide the number of publications by the total number of publications in the same year (recorded in Figure \ref{fig:glo}),
we end up with Figure \ref{fig:rcon}. It gives us more information about
the relative changes between countries.
For example the United States covers around 20\% of all publications,
but this ratio is slightly decreasing.
One clear increase is China, it also covers about 20\% of all publications since 2014
while this number was just 10\% in 2000.
Japan has also suffered from a clear relative decline between 2002 and 2006.
We may consider several natural questions related to this general picture,
for example what is the impact of the country of affiliation of the author(s)
on the citation of a paper\;\!? Again, this information has to be used
as a predictor in our analysis.

\subsection{The data}\label{sub:data}

As we mentioned in the Introduction, MathSciNet is not equipped with any tool
for extracting data from the website.
For that reason we had to carefully choose our data.
Since our initial motivation was to study the effect of international collaborations,
we fixed Japan as a target country. Also, in order to eliminate the time
variable from our investigations (this effect should be considered in the future),
we also fixed a publication year: 2009.
As a consequence, all the predictors are from 2009
while the numbers of citations have been collected at the end of 2018.

The total number of publications in 2009 with at least one author with a Japanese
affiliation is 3907 (in December 2018).
It is important to observe that around $55\%$ of these publications have zero or one
citation and that the number of works with many citations (for example more than 30)
is very limited. For that reason we stratify the data and draw a simple random sample
from each strata. Altogether, we shall use about $8\%$ of the total data
to do our analysis.  The information of the publications and stratified sample are
shown in Table \ref{tab:sam}.

\begin{table}
\begin{center}
\caption{Data collection}\label{tab:sam}
\begin{tabu} to 1.5\textwidth { | c | c | c | c | }
 \hline
Level & Number of citations & Number of publications & Number in sample \\
 \hline
5 & $>$ 30  & 54  & 4 \\
\hline
4 & 11 - 30 & 298 &  23 \\
\hline
3 &  6 - 10 & 413 & 23 \\
\hline
2 &  2 - 5 & 1008 & 91  \\
 \hline
1 &  1 & 628 & 48 \\
\hline
0 & 0 & 1506 & 116 \\
\hline
&Total& 3907 & 305 \\
\hline
\end{tabu}
\end{center}
\end{table}

\subsection{The predictors}\label{sub:pred}

Let us focus on the works published in 2009 with at least one author
affiliated with a Japanese institute. Note that for confidentiality reasons,
we collected only the countries of the institutes, or the numbers of different institutes,
but not the names of the institutes. Clearly, any information on the author(s)
is also disregarded.

The predictors are gathered in three families, namely: authors, article, journal.
In the subsequent lists, we have tried to use all information
provided by MathSciNet. Note that all the predictors are related to
the year of publication of the articles, namely 2009.

\begin{itemize}
\item Authors
\item[${}$] \emph{aut}: the number of authors
\item[${}$]\emph{rfstp}: proportion of young researchers (see below)
\item[${}$]\emph{avacag}: the average academic age, calculated by $\big[\sum\limits_{\hbox{authors}}(\hbox{2009 $-$ year of first publication})\big]$\Big/\emph{aut}
\item[${}$]\emph{inst}: the number of different institutes of authors
\item[${}$]\emph{nati}: the number of different nations of institutes
\item[${}$]\emph{jpper}: the proportion of authors from Japanese institutes
\end{itemize}

The proportion of young researchers corresponds to the ratio of the number of authors
who published their first work in 2009 divided by the total number of authors (which
is given by \emph{aut}).

\begin{itemize}
\item Article
\item[${}$]\emph{ref}: the number of references of the paper (as provided by MatSciNet,
with $0$ if no information is provided by MathSciNet)
\item[${}$]\emph{pg}: the number of pages of the article
\item[${}$]\emph{rt}: 1 (there is a review text in MathSciNet), 0
(there is no review text)
\item[${}$]\emph{msccla}: the combined subject classification used in Figure \ref{fig:pri}
\item[${}$]\emph{claper}: the proportion of each subject in \emph{msccla} during 2004 to 2008 computed by $$\frac{\hbox{the number of papers from 2004 to 2008 on subject }\emph{msccla}}{\hbox{the total number of papers from 2004 to 2008}}$$
\end{itemize}

Let us note that a review text is a summary of a paper written by a reviewer of
MathSciNet. Most of the works collected by MathSciNet are attached with a review text,
but some of them are not. If there exists a review text, we assign the number 1
to the predictor \emph{rt}, and we assign 0 to \emph{rt} if such a review text
does not exist.

\begin{itemize}
\item Journal
\item[${}$]\emph{mcq}: Mathematical Citation Quotient for 2009 (see below)
\item[${}$]\emph{cito2009}: the accumulated citations of the journal up to 2009
\end{itemize}

The Mathematical Citation Quotient (MCQ) for 2009 and for each journal
are provided by MathSciNet. For a given journal X  the MCQ for 2009 is
computed by the formula
\begin{align*}
\sum\limits_{i = 2004}^{2008}c_{i} \bigg/ \sum\limits_{i = 2004}^{2008}n_{i}
\end{align*}
where $c_{i}$ is the number of citations which appeared in 2009 about papers published
in year $i$ in the journal X, and $n_{i}$ is the number of publications in the journal X in
year $i$. For books or theses, \emph{mcq} and \emph{cito2009} are both 0.
In the sequel, we shall consider the predictor \emph{mcq}
as an approximate measure of the quality of a journal.
Since our aim is not to perform a comparison of the quality of the journals,
this rough information will be sufficient for our purpose.

For the response, we have collected the number of citations
provided by MathSciNet at the end of 2018.
These numbers have been distributed into 6 different levels,
according to the number of citations:
5 for more than 30 citations, 4 for a number of citations between 11 and 30,
3 for a number of citations between 6 and 10, 2 for a number of citations between
2 and 5, 1 for a number of citation equal to 1 and 0 for 0 citation.
These levels are reported in the first column of Table \ref{tab:sam}.

\section{Preliminary experiments}\label{pre}

In this section, we gather some preliminary results based on our sample.
The tree classifier will be introduced only in the next section.

\begin{wrapfigure}{R}{0.5\textwidth}
\centering
\includegraphics[scale =0.35]{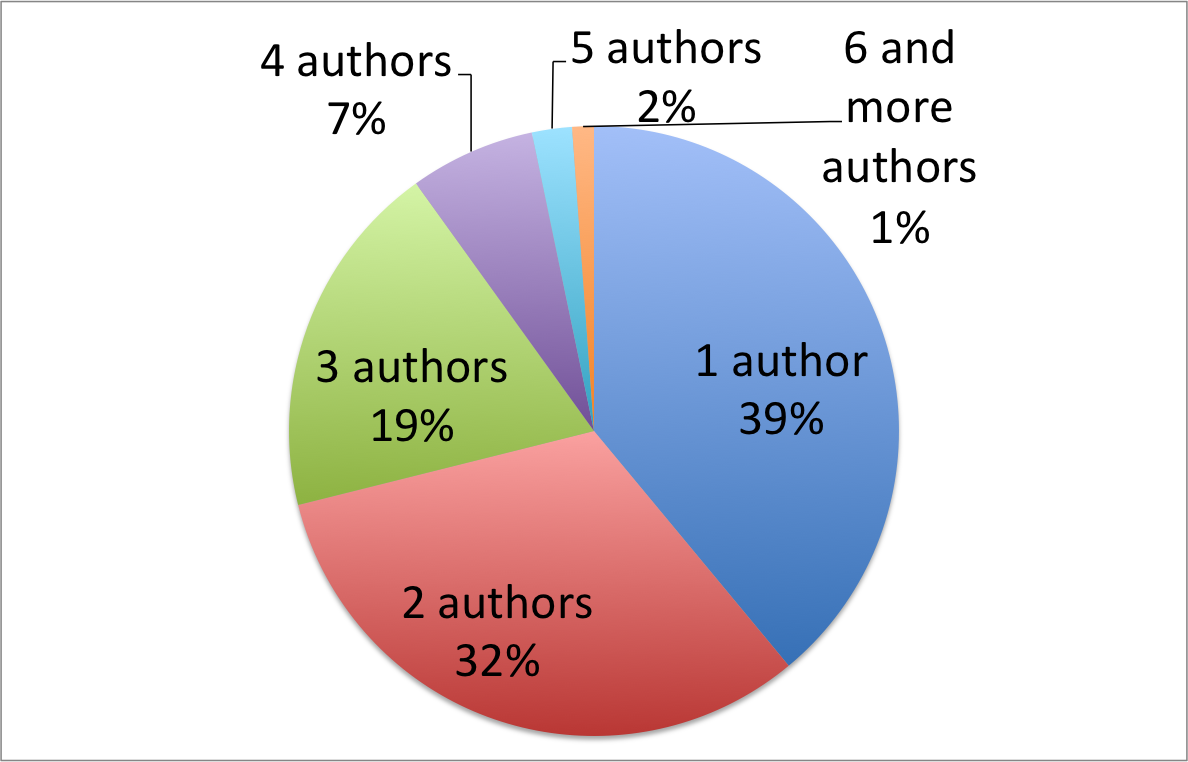}
\caption{Works produced in 2009 with at least one author from a Japanese institute}
\label{fig:jpaut}
\end{wrapfigure}

Let us recall that our initial motivation was whether collaborating with colleagues
from foreign institutes can increase the number of citations.
Since the works done by a single researcher are not collaborative works,
thus we should focus on the items written by two or more persons.
As shown in Figure \ref{fig:jpaut}, about $39\%$ of the 3907 papers is
made of publications written by only one author.
For an analysis about international collaboration versus no international collaboration
the corresponding items should not be used.
However, we do not want to do a census at this step, so let us keep all
data from our sample, and observe that the predictor \emph{aut} deals automatically
with the information on collaborative work or not.

\begin{figure}
\centering
\includegraphics[scale = 0.42]{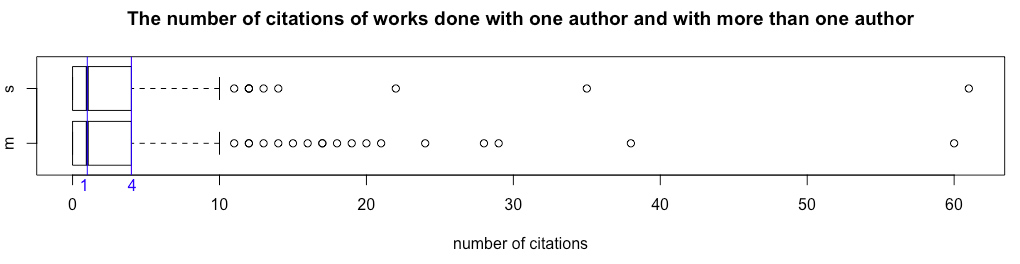}
 \caption{``s" means single author, ``m" means more than one author}
 \label{fig:bx3}
\end{figure}

Let us still investigate about the importance of one single author versus more than
one author. One could expect that a work written by a single author may have less
chance to be cited, compared to a collaborative work.
However, as the boxplot of Figure \ref{fig:bx3} shows,  the shape of the two categories
are almost the same: the minimum is zero, the median is one and the 3rd quartile is four. Based on this information, one may suspect that the number of authors
has a limited effect on the number of citations. Note that the graph with a single
author contains 102 observations, while the graph with more than one author contains 203 observations.

\begin{figure}
\centering
\includegraphics[scale = 0.41]{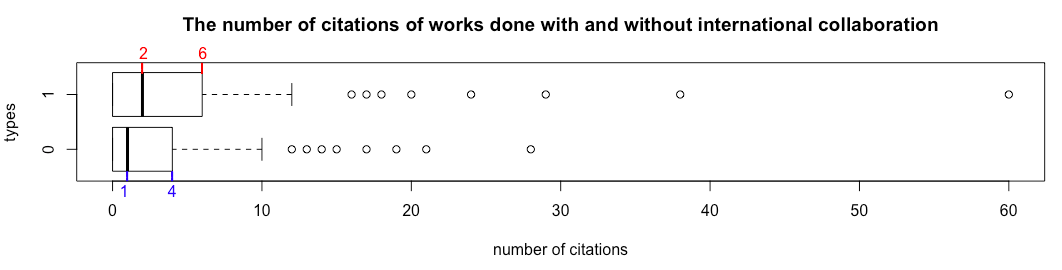}
 \caption{Type 0  =  local collaborations,  type 1 = international collaborations}
 \label{fig:bx1}
\end{figure}

If we now focus on international collaborations versus collaboration among
Japanese researchers, we can draw the boxplot of Figure \ref{fig:bx1}.
Type 0 contains all items with only Japanese authors (at least two), while type 1
contains items with at least one author from Japan and at least one author from abroad.
It is not a surprise that both categories have same minimum and 1st quartile which is
zero citation since more than one third observations have never been cited.
However, one can observe that the median of type 1 is 2 while the median of type 0
is 1. The 3rd quartile of type 1 is 6 which means $75\%$ of observations in type 1 have between zero and 6 citations.
For type 0, $75\%$ of items have between zero and 4 citations. Clearly type 0 is more concentrated while type 1 is more spread.
This is the first appearance of a distinction between local collaborations
and international collaborations. Note however that type 0 contains 138 observations
while type 1 contains only 65 observations.

\begin{figure}
\centering
\includegraphics[scale = 0.33]{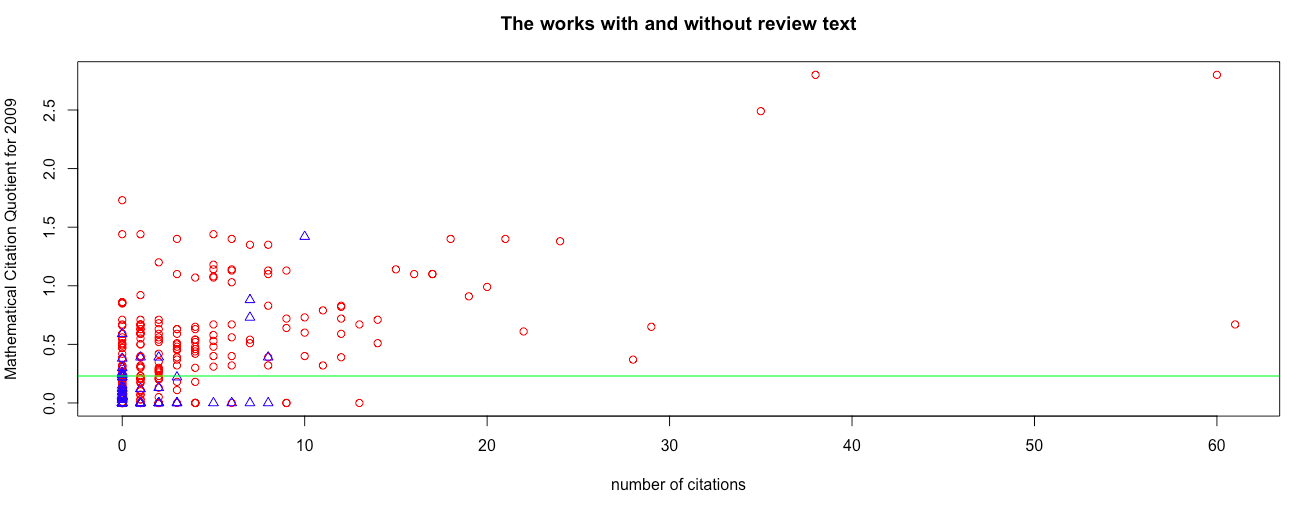}
 \caption{Triangle = absence of review text, circle = existence review text}
 \label{fig:bx4}
\end{figure}

During the data collection, we observed that MathSciNet contains quite a lot of items
without a review text (even though the idea of MathSciNet is that every work should be associated with a review text), and that some of these items were published on relatively less famous journals (low \emph{mcq}).
Based on this observation, one may want to know whether the existence
of a review text is linked to the quality of the journal, and also if it affects
the number of citations or not\;\!?  For answering this question we
realize the Figure \ref{fig:bx4}.
The horizontal line represents the average \emph{mcq} in 2009,
which is 0.23. We can see that most of the triangles
(the works without review text) are located under this line, which means these works
were published on relative low quality journals.
However, this last remark does not apply to items with
\emph{mcq} equal to 0 since these publications often appear in
books or in proceedings of conferences, for which the \emph{mcq} is 0
as explained in Section \ref{sub:pred}.

For the circles (the works with review text), though the number of citations has a
vague positive correlation to the \emph{mcq}, there are several items which were
published on high \emph{mcq} journals, and still have zero citation.
This observation means that the good quality of a journal can not promise numerous citations for a given paper appearing in this journal.

\section{Investigations with a classification tree}\label{invest}

Based on the preliminary observations performed so far, we have a better
feeling about the data. Clearly, our data set consists in 13 predictors, with
a mixture of quantitative and qualitative features. Also, some features that
were considered to be critical may be less important than what we thought,
and for sure there is no linear relations between the predictors and the response.
To deal with the special structure of this sample, we will use the tree-based
method which is considered to be a reasonable choice for analyzing this
kind of data.

\subsection{The tree classifier}

For this approach, let us just recall its main idea.
Suppose that  $\X$ be a $p$-dimensional measurement space
which contains all predictor vectors, and let $x = (x_1, x_2, \dots, x_p)$ denote
the elements of $\X$. Let $\C$ denote the set of possible responses, and for
simplicity assume that $\C$ is a finite set: $\C:=\{1,2,\dots,J\}$.
Each element of $\C$ is called a class.
A \emph{classifier} or \emph{classification rule} is a function
$d : \mathcal{X} \ni x\mapsto d(x)\in \C$.

\begin{figure}
\centering
\includegraphics[scale = 0.36]{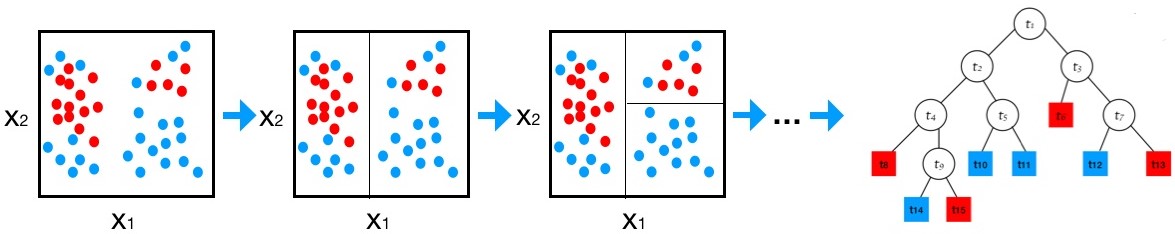}
\caption{Tree-structured classifier}
\label{fig:tr}
\end{figure}

If we suppose that $\X$ is only $2$-dimensional and if we also assume that
$\C$ contains only two elements, which are represented
by red (dark) and blue (light) dots in Figure \ref{fig:tr}, then the tree-based method consists
in separating the observations into a finite number of subregions.
The split points correspond to some values of the predictors,
and each successive subdivision is chosen according to a certain rule.
More precisely, each split needs to result in two subregions with lower impurity,
the impurity function is defined in Definition \ref{def:imp}.
The procedure will firstly run over all possible splits of one predictor,
then go to the next predictor, run over all the possible splits,
then go through all predictors.
It will select the best split of the best predictor, which means choose the predictor
and the split which lead to the biggest decrease of the impurity.
For each subregion, the classifier can keep separating them until meeting some
stopping rule.

\begin{definition}\label{def:imp}
An \emph{impurity function} is a function $\phi$ defined on all J-tuples
$(p_{1}, \dots, p_{J})$, satisfying $p_{j} \geq 0$ for $j = 1, \dots, J$  and
$\sum_{j} p_{j} =1$, with the following properties:
\begin{itemize}
\item[$(1)$] $\phi$ is a maximum only at the point $\Big(\frac{1}{J}, \frac{1}{J}, \dots, \frac{1}{J}\Big)$,
\item[$(2)$] $\phi$ achieves its minimum only at the points $(1, 0, \dots, 0)$, $(0, 1, 0, \dots, 0)$, $\dots$, $(0, \dots, 0, 1)$,
\item[$(3)$] $\phi$ is a symmetric function of $p_{1}, \dots, p_{J}$.
\end{itemize}
\end{definition}

When the procedure stops, we have a tree-structured classifier. Each circle or square
represented in Figure  \ref{fig:tr} is called a node, and the terminal nodes are called leaves.

\subsection{Classification tree}\label{sec:class}

In this section, we use R and the tree-based methods to analyze the sample given in Section \ref{sub:data}.

In Table \ref{tab:sam}, one can see that despite the papers published in 2009
have been distributed in 6 levels only, the six populations are not equally distributed.
For example, there are only 54 items with more than 30 citations, which represents
less than $2\%$ of all the works in 2009. On the other hand,
nearly $40\%$ of all publications are at level 0 which corresponds to no citation.
It is known that accurately detecting \emph{rare events} is one difficulty
when we construct the classifier. Firstly, because a limited number of examples may
not be enough to train the classifier, and secondly, because the information
contained in the sample may not be enough to distinguish rare events from
more common events. To deal with this situation, several strategies have been developed, for example the oversampling and undersampling inspired from signal processing, or the SMOTE \cite{CFSA} which consists in artificially creating synthetic data points to balance the distribution of the sample, or collect more data in these classes.
In these preliminary investigations, we shall not go further in these directions.

During our investigations, we also observed that the procedure does not easily distinguish
items of level 0 with items of level 1. One natural explanation is that works
which are cited 0 or 1 time have quite similar predictor vectors.

Based on these observations, and since our sample is not so big,
we decided to combine some levels together, and ended up
with a response consisting of 3 classes only:
the \emph{low level} with 0 or 1 citation,
the \emph{median level} with a number of citations between 2 and 10,
and the \emph{high level} corresponding to papers with 11 or more citations.
On the other hand, for the predictors we have implemented the 13 predictors
introduced in Section \ref{sub:pred}.

In R, the tree classification code is {\fontfamily{qcr}\selectfont tree}.
Then, one needs to use the command {\fontfamily{qcr}\selectfont control}
to control the growing of the tree classifier. The {\fontfamily{qcr}\selectfont nobs}
is the number of observations in the training set, which is equal to 305 in our case.
The {\fontfamily{qcr}\selectfont mincut} controls the minimum number of observations
to include in each child node. Namely, if a split can not result in both nodes with
more than {\fontfamily{qcr}\selectfont mincut} observations, then the split would
not be carried out. Similarly, the {\fontfamily{qcr}\selectfont minsize} is the minimum number of observations in any parent node.  In our case, to avoid reaching the
maximum depth of growing the tree, we choose 3 and 6 respectively.
The splitting rule is implemented by choosing the split point which minimize
the \emph{deviance}.

The deviance of each node $t$ (including the terminal nodes), denoted by $D_{t}$, is computed by the following formula:
\begin{equation*}
-2\sum_{j=1}^J p_{j \mid t}n_{t}\ln p_{j \mid t},
\end{equation*}
where $n_{t}$ is the number of observations in node $t$, and $p_{j \mid t}$ is the proportion of observations with class $j$ in node $t$.
Note that if we set $\phi(p_{1\mid t}, \dots, p_{J \mid t}) = D_{t}$, then
$\phi$ is concave and meets the three properties in Definition \ref{def:imp}.

Let us also define \emph{residual mean deviance} by the following formula:
\begin{align}\label{rmd}
\frac{\sum\limits_{t \in T'}D_{t}}{N - |T'|},
\end{align}
where $D_{t}$ is the deviance of terminal node $t$,
$T'$ denotes the set of terminal nodes (leaves) of a tree $T$,
$|T'|$ corresponds to the number of leaves of the tree $T$, and
$N$ is the total number of observations.
This residual deviance is a global property of the tree-structured classifier,
and a small number is preferred (among all possible classifier trees).
As a consequence, we always look for a small numerator in \eqref{rmd},
and concomitantly we wish to get a tree with simple structure.
Indeed, such a tree has less terminal nodes, and therefore a bigger denominator in
\eqref{rmd}.

After constructing the tree-structured classifier based on our sample,
R reports a summary of this classification as follows:

{\fontfamily{qcr}\selectfont
\noindent
    >Variables actually used in tree construction:\\
        "ref"      "mcq"      "pg"       "rfstp"    "msccla1"  "avacag"   "rt"      \\
        "aut"      "nati"     "cito2009"\\
     >Number of terminal nodes:  62 \\
     >Residual mean deviance:  0.6544 = 159 / 243 \\
     >Misclassification error rate: 0.1311 = 40 / 305 \\
}

This summary tells us that for the construction of the tree, only 10 predictors
have been used. It also provides the residual mean deviance, and the misclassification
error rate.
The misclassification error rate is the ratio of the number of incorrectly classified items and the total number of items in the sample. Indeed, since the tree assigned to each leaf
a class low level, median level, or high level, the misclassified items in each leaf can
easily be identified. From the summary we see that 40 observations are misclassified, which is about $13\%$ of the sample. Table \ref{tab:sum} gives the details of classification, the colored cells represent the correct classifications. One should notice that the misclassification rates for high level and median level items are higher than $13\%$.

\begin{table}
\begin{center}
\caption{The summary of classification}\label{tab:sum}
\begin{tabu} to 1.5 \textwidth { | c | c | c | c | c | }
 \hline
 \diagbox[width=10em]{predict value}{real response} & high level & median level & low level &  predict error rate \\
 \hline
 high level  & \cellcolor{magenta!60}21  & 2 & 3 & 0.19 \\
\hline
 median level & 4 &  \cellcolor{magenta!60}92 & 10 & 0.13 \\
\hline
 low level & 2 & 20 & \cellcolor{magenta!60}151 & 0.13 \\
\hline
 misclassification rate & 0.22 & 0.19 & 0.08 & N/A \\
\hline
\end{tabu}
\end{center}
\end{table}

Note that among the 62 terminal nodes, 30 are pure (which means
that all items are of the same class) and these pure leaves perfectly classify 175 observations in the dataset. In summary, the misclassification rate reported by
R may not be evenly spread in the feature space, to look into the distribution of
each node may help us better understand the performance of the procedure.
Also, if one can detect the region in which the classifier performs improperly,
then finding and interpreting the reason may be very useful for the further investigation.

As mentioned above, not all predictors are used in the construction of the tree: \emph{claper}, \emph{jpper} and \emph{inst} are not selected as useful predictors.
Recall that \emph{claper} is the proportion of each subject in \emph{msccla}.
As a consequence, we interpret this result as the fact that that \emph{claper}
can not make any efficient distinction in any part of the sample compared
with other predictors. Similar interpretation can be given for \emph{jpper} and
\emph{inst}, which means one can also not predict the number of citations by
considering the proportion of researchers from Japanese institutes and from
the number of distinct institutes.

It is natural to say that the abandoned predictors are not as important as others.
On the other hand, how should we order the selected predictors by their relative
importance, and how should we define the importance of a predictor\;\!?
For example, if one looks at each split point in the tree classifier,
one finds that the number of times each predictor is used at a split point are
quite different. Such numbers are reported in Table \ref{tab:split}.
Does it mean an often used predictor is more important than others\;\!?
What about the decrease of the deviance at each split\;\!? Should we also consider the ability of decreasing the deviance\;\!? We will discuss these in the next section.

\subsection{Importance of predictors}\label{sec:imp}

Our aim in this section is to discuss the relative importance of each predictor.
Several approaches for answering this question are possible, and we shall discuss
two of them.

Firstly, let us observe that the number of times a predictor defines a split
point is not directly related to its importance. This may be affected by the nature
of the predictors, as for instance the review text \emph{rt}. There are only two
values in \emph{rt}, which are zero and one. If the predictor was used at the
very beginning, then it would never be used again because one can not make any
further division on the space of the predictors.

Secondly, the growth of the tree is controlled by the decrease of the deviance.
The tree aims to divide the sample into smaller parts with lower deviance.
One natural idea is to consider the ability of decreasing the deviance,
the one which decreases the deviance most would be the important one.
However, like we said in the last paragraph, this may be affected by the intrinsic
properties of the predictors. It would be more convincing to use the average
decrease of each split.

In Table \ref{tab:split} the third column gives the information on the total decrease
of deviance due to each predictor, and the forth column gives the average
decrease for each split. The ranking of importance of predictors with respect to the
average decrease of deviance is shown in the fifth column called rank (dev).

\begin{table}
\begin{center}
\caption{Importance of predictors: deviance drop versus misclassification rate}\label{tab:split}
\begin{tabu} to 1.5\textwidth { | c | c | c | c | c || c | c |}
\hline
predictors & \# of splits & Deviance drop & Deviance drop $/$ & rank & increase & rank \\
 & & & per split  & (dev) & misc. rate  &  (mis)\\
\hline
\emph{avacag}  & 15  & 72.7 & 4.8 & 7 & 0.145  & 3\\
\hline
\emph{pg} & 12 & 63.2  & 5.3 & 6 & 0.120 & 4\\
\hline
\emph{mcq} & 8 & 53.5  & 6.7 & 3 & 0.201 & 2\\
\hline
\emph{msccla} & 7 & 40.9  & 5.8 & 4 & 0.095 & 5\\
\hline
\emph{aut} & 6 & 20.5  & 3.4 & 9 & 0.027 & 8\\
\hline
\emph{cito2009} & 5 & 43.2 & 8.6 & 2 & 0.077 & 6\\
\hline
\emph{ref} & 4 & 89.3  & 22.3 & 1 & 0.211 & 1\\
\hline
\emph{rfstp} & 2 & 11.3 & 5.7 & 5 & 0.032 & 7\\
\hline
\emph{rt} & 1 & 3.9  & 3.9 & 8 & 0.007 & 9\\
\hline
\emph{nati} & 1 & 1.2  & 1.2 & 10 & 0.004 & 10\\
\hline
Total & 61 & 399.7 & 6.6 & &  & \\
\hline
\end{tabu}
\end{center}
\end{table}

Another direction of considering this question is to put incorrect information
separately for each predictor and record the change of misclassification rates.
In general, the incorrect information will result in worse classification accuracy.
In this research, we will use a technique called \emph{shuffle} to provide the
importance of predictors.
The idea of shuffle is very simple. After growing the tree, the classification rules
are fixed and the misclassification error rate is $\beta$.
Now, we choose one predictor, say \emph{ref}, and randomly disorder the value
of \emph{ref} for each observation in the sample, then we obtain an artificial test
dataset. Let us use the tree for classifying test data, and record the misclassification
rate (it means the total rate of error in the leaves). We call this misclassification rate
$\beta_{k}$, and repeat these steps for $K$ times. The average increase of misclassification
rate $\gamma_{\emph{ref}}$ is given by
\begin{align*}
\gamma_{ref} := \frac{1}{K}\sum_{k=1}^K(\beta_{k} - \beta).
\end{align*}

Next, we move to the next predictor and repeat the above-mentioned steps.
Finally, we record these average increases of misclassification rate for each predictor
in the sixth column of  Table \ref{tab:split}.
For our experiment we have fixed $K=1000$.
As a final step, we rank the predictors according to the increase
of the misclassification rate, as shown in the seventh column
of Table \ref{tab:split}, called rank (mis).
The idea behind this classification is that random data for a non-important predictor
should not change the classification error rate a lot, while any randomness for an important
predictor would increase much more the misclassification rate.

The predictors \emph{ref} and \emph{mcq} cause more than $20\%$ of increase of misclassification rate on average. The predictors \emph{aut}, \emph{rt}
and \emph{nati} are still the last three predictors: disordering their values cause only
a very limited increase in the misclassification rate.

Compared with the previous ranking (the ability of decreasing the deviance) one
observes that the new ranking is rather different, except for its extremes (number 1
and number 10) . For example the ranking of the average academic age \emph{avacag}
and of the number of pages \emph{pg} are higher with the new method.
This can partially be explained by looking again at Table \ref{tab:split}:
\emph{avacag} and \emph{pg} are the most frequently used splitting predictors,
which means that incorrect information of \emph{avacag} and \emph{pg} will
generate some misclassifications several times. In other words, the wrong
information will be checked many times and each time the procedure will move
away from the correct classification.

\section{Conclusion}

The research was inspired by a rumor, but the initial question was very quickly
extended to a more general analysis of the predictors related to the publication
of mathematics papers. Data have been extracted from MathSciNet,
but the limitation with the current access to the database has been major limitation
for our investigations.
Nevertheless, various information have been computed from the raw data,
and several figures have been created for illustrating the outcomes.
Because of the diversity of the predictors, a tree-based
method has been chosen, and the predictors have been ranked according
to two different criteria.

Two of the rather surprising outcomes from our investigations and based on
our sample are : 1) the weak impact on the number of authors (\emph{aut})
for the number of citations,
2) the relatively weak impact of the number of nations involved
in a research project (\emph{nati}) again for the number of citations.
These results do not support the conclusion obtained in \cite{SWBA}
based on a multidisciplinary approach.
In that respect, further investigations are certainly necessary.
Note however that our results can not be directly compared with this
reference since we are taking a much larger set of predictors into account.
With more data, it would also be necessary to consider more classes for the response,
and therefore get a finer analysis of the number of citations.
For example, the small difference seen in Figure \ref{fig:bx1} about the
number of citations between the works done with and without international collaborations
is completely invisible in the analysis with the three classes since
this difference vanishes in the single \emph{median level}.

It is also interesting to recall that three predictors have not been
used for the elaboration of the tree, namely the one related to
the number of different institutes of authors (\emph{inst}),
the one related to the proportion of authors from Japanese institutes
(\emph{jpper}), and the one computed from the combined
Mathematics Subject Classification (\emph{claper}).
Even so the lack of importance of \emph{jpper} is natural,
the absence of \emph{claper} in the tree might be due
to our combined version of the subject classification.
Further investigations with the full power of the MSC would certainly be
instructive.

In Section \ref{sec:imp}, we provided two rankings of the importance of predictors,
the first one is based on the average deviance drop, and the second one is defined
by using a shuffle technique. The results show that the number of references
\emph{ref} and the Mathematical Citation Quotient \emph{mcq} are important
predictors. For the second one, it is not surprising since the \emph{mcq} is a measure
of the quality of a journal. Nevertheless, this information is important on average
for papers published in a journal, but does not always mean anything for a precise
paper, as emphasized by Figure \ref{fig:bx4}.
The importance of \emph{ref} is also not so surprising (and had already been observed
by other researchers like in \cite{SWBA}), since many references might
either reflect the importance of a new contribution, or will be more easily visible
by a large number of cited researchers. Note that the predictor
\emph{cito2009} providing the accumulated citations of a journal up to 2009
also reflects the importance of a journal, but is more difficult to interpret due to
the difference of age between journals.

For the tree classifier, the selection of impurity function affects the growth of the tree.
A good impurity function can generate descendants with more pure nodes and thus less efforts are necessary for the subsequent splits.  In this research, the deviance
has been chosen as the splitting rule, and it results in a misclassification error rate
of $13\%$ with 62 terminal nodes.
However, the size of the tree is certainly too big compared with the size of the sample,
and once again we need more data for a more precise analysis.
Let us also mention that the performance of classification on low and median levels
items is better than it on high level items. We can say that the classification abilities
of the classifier on different parts of the feature space are different.
Some parts of the sample are easily classified, but other parts are hard to classify
for many reasons. Finding methods to increase the classification accuracy in those
parts should be considered in a future study. One possible solution is to take more 
data such that it can provide more information to the model to learn the full picture
of the dataset. On the other hand, advanced algorithms may improve the accuracy and stability of the tree classifier. For example, Mahmood introduced two algorithms
aimed to improve the accuracy in \cite{MRR}, and Zhang and Jiang presented a
splitting criteria based on similarity in \cite{ZHJ}. For increasing the stability, two famous improvements done by Breiman in \cite{B1,B2} are bagging predictors and random forests.
In \cite{LMM}, an information-theoretic method was introduced for building stable and comprehensible decision tree models.

Let us finally briefly mention some additional investigations which have
been performed but not reported in this paper.
We also tried to look at the performance of the classifier on a test data.
For that purpose, one could have randomly divided the sample into training
set and test set, then used the training set to build the tree and let  the test
set go through the tree, and record the test error.
This error is computed from the number of misclassified items divided by the total
number of items in the test set.
In general, the test error will be greater than the misclassification error rate of the tree built
with all data in the sample. Since we use less data, the classifier can not grab all the information in the sample. And this situation would be worse when the sample
size is small. Also, since the division of the sample is random,
the result will vary depending on the division results. For this reason, we usually
divide the sample several times, and build many trees using the training sets.
We then put the test sets into the corresponding trees respectively and record the test errors. We finally report the average test error by using each test error in the pre-mentioned procedures.
For our experiment, we used 1 item as test set and the left 304 items as training set,
we run the test for 1000 times and the average test error was $45.4\%$. Since the response consists of three levels, when there is no prior information of the distribution of the three levels, take a random guess would give you an error rate of $66.6\%$. The tree classifier gives $21\%$ increase of classification accuracy.

In mathematics, the motivating rumor is certainly not as strong as in other
disciplines.
Except Figure \ref{fig:bx1}, there is no significant evidence showing that working with researchers from foreign institutes will be helpful for increasing the number of citations.
However, our investigations have suffered from the strict rules of access to the
database, and our preliminary results need further investigations and confirmations.
There is still a lot to explore in MathSciNet, and we hope to do it in the near future.

\end{document}